\documentclass[smallcondensed]{svjour3-id}

\usepackage{latexsym, amsfonts, amsmath, amssymb, verbatim}
\usepackage{tikz, float}
\usepackage[top=3cm, bottom=3cm, left=3cm, right=3cm]{geometry}
\usetikzlibrary{plotmarks}

\newcommand{\AC}{{AC}}
\newcommand{\BV}{{BV}}

\newcommand{\CTPP}{\mathrm{CTPP}}

\newcommand{\Lip}{{\mathrm{Lip}}}

\newcommand{\mR}{\mathbb{R}}
\newcommand{\mC}{\mathbb{C}}

\newcommand{\abs}[1]{\left\lvert#1\right\rvert}

\newcommand{\vecx}{{\boldsymbol{x}}}
\newcommand{\vecy}{{\boldsymbol{y}}}
\newcommand{\vecz}{{\boldsymbol{z}}}

\newcommand{\vecq}{\boldsymbol{q}}

\newcommand{\vecu}{{\boldsymbol{u}}}

\newcommand{\vecw}{{\boldsymbol{w}}}

\newcommand{\norm}[1]{\left\lVert#1\right\rVert}

\newcommand{\normbv}[1]{\left\lVert#1\right\rVert_{\BV(\sigma)}}

\renewcommand{\Re}{\mathop{\mathrm{Re}}}
\renewcommand{\Im}{\mathop{\mathrm{Im}}}

\newcommand\st{\thinspace : \thinspace}
\def\ls[#1,#2]{\overline{\vphantom{\vbox to 1.2 ex{}} #1\, #2}}

\def\sparen(#1){\Bigl ( #1 \Bigr )}
\def\ssparen(#1){ (#1) }
\newcommand\plist[1]{\bigl[ #1 \bigr]}

\DeclareMathOperator*{\var}{var}
\DeclareMathOperator{\cvar}{\rm cvar}
\DeclareMathOperator{\vf}{vf}

 \newcommand{\mysubsection}[1]{\subsection{#1}\ \smallskip\noindent}

\begin{document}

\title{Approximation in $AC(\sigma)$}  
\author{Ian Doust \and Michael Leinert \and Alan Stoneham}

\institute{Ian Doust \and Alan Stoneham
\at School of Mathematics and Statistics,
University of New South Wales,
UNSW Sydney 2052 Australia \\ 
\email{i.doust@unsw.edu.au}, \email{a.stoneham@unsw.edu.au}
\and
 Michael Leinert
 \at  Institut f{\"u}r Angewandte Mathematik,
Universit{\"a}t Heidelberg,
Im Neuenheimer Feld 205,
D-69120 Heidelberg
Germany \\ 
\email{leinert@math.uni-heidelberg.de}}
\date{}


\maketitle

\begin{abstract}
For a nonempty compact subset $\sigma$ in the plane,
the space $AC(\sigma)$ is the closure of the space of complex polynomials in two real variables under a particular variation norm. In the classical setting, $AC[0,1]$ contains several other useful dense subsets, such as continuous piecewise linear functions, $C^1$ functions and Lipschitz functions. In this paper we examine analogues of these results in this more general setting.
\end{abstract}

\section{Introduction}
A standard method in analysis is to first prove some fact for particularly well-behaved functions, and then extend the result to a more general class via a limiting argument. In order to do this one  must of course know that the functions that you are interested in can be approximated by the simpler well-behaved ones.

The class of $AC(\sigma)$ spaces, comprising `absolutely continuous' functions   on a nonempty compact subset $\sigma$ of the plane, was introduced in \cite{Ash} and \cite{AD1} to provide a functional calculus model which unified the theories for the well-bounded operators of Smart and Ringrose \cite{Sm,Rin1,Rin2} and the trigonometrically well-bounded operators of Berkson and Gillespie \cite{BG2}, while overcoming some of the problems inherent in the definition of $AC$-operators \cite{BG1,BDG}. Much is now known about the algebraic structure of these spaces \cite{AS1,AS2,AS3,DA,DL}, but rather less has been recorded about the relationship between these spaces and other standard spaces of functions.

Throughout the paper we will freely swap between considering the plane as $\mC$ or $\mR^2$. Subsets of $\mR$ may be considered as embedded in $\mC$ or as subsets of the first coordinate axis in $\mR^2$ as appropriate. Given a nonempty subset $\sigma \subseteq \mR^2$, the space $AC(\sigma)$ is defined to be the closure of the complex polynomials $p(x,y)$ in two real variables inside the space $BV(\sigma)$ of functions of bounded variation (in the sense of \cite{AD1}).

In the classical situation of $AC[0,1]$, it is useful to know that this space contains several standard spaces as dense subsets; $C^1[0,1]$, the Lipschitz functions, and the continuous piecewise linear functions. The aim of this paper is to examine to what extent the analogous facts are true in this more general setting.

In Section~\ref{S:def} we shall briefly give the main definitions of the spaces involved. A more detailed account can be found in \cite{DLS}. Note that the definitions given here involve a significant simplification of the original ones used in \cite{AD1}. A discussion of the equivalence of the different definitions can be found in the appendix to \cite{AS3}.

In Section~\ref{S:CTPP} we shall introduce the space $\CTPP(\sigma)$ of continuous piecewise planar functions on $\sigma$, which forms a natural analogue of the continuous piecewise linear functions on the line.
An important fact which we shall use in this section is that being absolutely continuous is a `local property' in  the sense that a function is absolutely continuous on $\sigma$ if and only if it is absolutely continuous on some neighbourhood of each point in $\sigma$.

In the final section we shall show that (suitably interpreted), we always have the inclusions $C^1(\sigma) \subseteq AC(\sigma) \subseteq C(\sigma)$. (These inclusions are of course in general proper, since that is the case when $\sigma = [0,1]$.)

\section{Definitions}\label{S:def}

Suppose that $\sigma$ is a nonempty subset of the plane and that $f: \sigma \to \mC$. The main concept of interest to us will be the variation of $f$ over $\sigma$. The definition looks reminiscent of the usual definition of variation for a function defined on an interval, but must now take account of the fact that there is no longer a natural ordering to the points of the set.

Suppose that $S =
\plist{\vecx_0,\vecx_1,\dots,\vecx_n} = [\vecx_j]_{j=0}^n$ is a finite ordered list of points in the plane.
Note that the elements of such a list do not need to be distinct. For the moment assume that $S$ contains at least two distinct points.
Let
$\gamma_S$ denote the piecewise linear curve joining the points of $S$ in order.

\begin{definition}\label{defn:x-seg}
Suppose that $\ell$ is a line in the plane. We say that $\ls[\vecx_j,\vecx_{j+1}]$, the line segment joining $\vecx_j$ to $\vecx_{j+1}$, is a \textbf{crossing segment} of $S$ on $\ell$ if any one of the following hold:
\begin{enumerate}
  \item $\vecx_j$ and $\vecx_{j+1}$ lie on (strictly) opposite sides of $\ell$;
  \item $j=0$ and $\vecx_j \in \ell$;
  \item $j > 0$, $\vecx_j \notin \ell$ and $\vecx_{j+1} \in \ell$.
\end{enumerate}
\end{definition}


Let $\vf(S,\ell)$ denote the number of crossing segments of $S$ on $\ell$ and define the \textbf{variation factor} of $S$ to be
 \[ \vf(S) = \max_{\ell} \vf(S,\ell). \]
To cover the case where $S$ contains a single point we set $\vf([\vecx_0],\ell) = 1$ if $\vecx_0 \in \ell$ and zero otherwise, and set $\vf([\vecx_0]) = 1$.
Note that $\vf(S)$ is well-defined since $1 \le \vf(S) \le n$. Informally, $\vf(S)$ may be thought of as the maximum number of times any line crosses $\gamma_S$.

It is clear that the variation factor is invariant under rotations or translations of the set $S$, and indeed under any affine transformation of $S$.

We define the \textbf{curve variation of $f$ on the list $S$} to be
\begin{equation} \label{lbl:298}
    \cvar(f, S) =  \sum_{i=1}^{n} \abs{f(\vecx_{i}) - f(\vecx_{i-1})}.
\end{equation}
We include the case $S = \plist{\vecx_0}$ by setting $\cvar(f, \plist{\vecx_0}) = 0$. The \textbf{two-dimensional variation} of a function $f : \sigma
\rightarrow \mC$ is defined to be
\begin{equation} \label{lbl:994}
    \var(f, \sigma) = \sup_{S}
        \frac{ \cvar(f, S)}{\vf(S)},
\end{equation}
where the supremum is taken over all finite ordered lists of elements of $\sigma$.
The \textbf{variation norm} is
  \[ \normbv{f} = \norm{f}_\infty + \var(f,\sigma) \]
and this is used to define the set of functions of bounded variation on $\sigma$,
  \[ \BV(\sigma) = \{ f: \sigma \to \mC \st \normbv{f} < \infty\}. \]

We record for later the most important properties of these spaces. Proofs can be found in \cite{DLS}.
\begin{enumerate}
\item $BV(\sigma)$ is a Banach algebra.
\item\label{consistency} If $\sigma = [a,b]$ is an interval, then $\normbv{f}$ is just the usual variation norm.
\item\label{real-imag} $f \in BV(\sigma)$ if and only if $\Re(f)$, $\Im(f) \in BV(\sigma)$.
\item\label{aff-inv} (Affine invariance) If $\phi$ is an invertible affine transformation of the plane and $\sigma' = \phi(\sigma)$, then $\Phi(f) = f \circ \phi^{-1}$ is an isometric isomorphism from $BV(\sigma)$ to $BV(\sigma')$.
\item\label{const-lines} Let $\sigma_\mR = \{\Re z \st z \in \sigma\}$ and suppose that $g \in BV(\sigma_\mR)$. Then the function $f(z) = g(\Re z)$ is in $BV(\sigma)$ and $\normbv{f} \le \norm{g}_{BV(\sigma_\mR)}$. By affine invariance, this result can be applied to any function which is constant on a family of parallel lines.
\item\label{poly-in-bv} $BV(\sigma)$ always contains the polynomials in two real variables (considering $\sigma$ as a subset of $\mR^2$).
\end{enumerate}

It follows from (\ref{poly-in-bv}) that the closure of the polynomials in two real variables is always a closed subalgebra of $BV(\sigma)$. This subalgebra is called the algebra of absolutely continuous functions on $\sigma$ and is denoted $AC(\sigma)$. Since convergence in the $BV$ norm implies uniform convergence, it is clear that every absolutely continuous function is in fact continuous. Again, in the case that $\sigma = [0,1]$, the space defined here is the same as the usual space $AC[0,1]$. (This is indeed true for  any compact subset of $\mR$. The spaces $AC(\sigma)$ for $\sigma \subseteq \mR$ were introduced by Saks \cite{S}, if not earlier.) Properties (\ref{real-imag}), (\ref{aff-inv}) and (\ref{const-lines}) also hold for $AC(\sigma)$ spaces.

It is easy to construct examples where a function is of bounded variation on each of two compact sets, but not of bounded variation on their union.  The following result from \cite{AD3} gives a sufficient condition on the subsets for the `joined function' to be of bounded variation.

\begin{theorem}[{\cite[Theorem 3.1]{AD3}}]\label{join-convexly}
Suppose that $\sigma_1,\sigma_2 \subseteq \mC$ are nonempty compact sets which are disjoint except at their boundaries. Suppose that $\sigma = \sigma_1 \cup \sigma_2$ is convex and that $f: \sigma \to \mC$. If $f|\sigma_1 \in BV(\sigma_1)$ and $f|\sigma_2 \in BV(\sigma_2)$ then $f \in BV(\sigma)$ with
  \[ \normbv{f} \le \norm{f|\sigma_1}_{BV(\sigma_1)} +
                       \norm{f|\sigma_2}_{BV(\sigma_2)}. \]
\end{theorem}

An important fact is the following result from \cite{DLS} that says that being absolutely continuous is a local property. We shall say that a set $U$ is a \textit{compact neighbourhood} of a point $\vecx \in \sigma$ (with respect to $\sigma$) if there exists an open neighbourhood $V\subseteq\mC$ of $x$ such that $U = \sigma \cap \overline{V}$.

\begin{theorem}[Patching Lemma]\label{patching-lemma}
	Suppose that $f : \sigma \to \mC$. Then $f \in \AC(\sigma)$ if and only if for every point $\vecx \in \sigma$ there exists a compact neighbourhood $U_\vecx$ of $\vecx$ in $\sigma$ such that $f|U_\vecx \in \AC(U_\vecx)$.
\end{theorem}

We shall denote the Lipschitz functions on $\sigma$ by $\Lip(\sigma)$. This is a Banach algebra under the norm $\norm{f}_{\Lip(\sigma)} = \norm{f}_\infty + L_\sigma(f)$, where
  \[L_\sigma(f) =
      \sup \left\{ \frac{|f(\vecx) - f(\vecx')|}{|\vecx-\vecx'|}
                 \st \vecx \ne \vecx' \in \sigma \right\}. \]
(For the degenerate case of a singleton set, we set $L_\sigma(f)$ to be zero for all $f$.)
In the following proposition \cite[Theorem 3.11]{DLS}, $C_\sigma$ stands for the variation constant of $\sigma$, which is the variation of the identity function $\zeta(\vecz) = \vecz$ over $\sigma$ (considered as a subset of $\mC$). If $\sigma$ lies inside a rectangle, then $C_\sigma$ is at most the sum of the width plus the height of the rectangle.

\begin{proposition}\label{Lip-functs}
If $f \in \Lip(\sigma)$ then $\var(f,\sigma) \le C_\sigma L_\sigma(f)$ and so
$\normbv{f} \le \max \{1,C_\sigma\} \norm{f}_{\Lip(\sigma)}$. Thus $\Lip(\sigma) \subseteq BV(\sigma)$.
\end{proposition}

Unfortunately, depending on $\sigma$, there may be Lipschitz functions which are not absolutely continuous. Examples are given in \cite{AD1} and \cite{DLS}.

Continuous piecewise linear functions on an interval are of course Lipschitz and hence absolutely continuous. The analogue for subsets of the plane are continuous piecewise planar functions.  We shall denote the space of continuous piecewise planar functions on $\sigma$ by $\CTPP(\sigma)$. The main part of this paper is to show that this space is always a dense subset of $AC(\sigma)$.

For $k  = 1,2,3, \dots$, let $C^k(\sigma)$ denote the algebra of functions $f$ for which there is an open neighbourhood $U$ of $\sigma$ and a $k$-times continuously differentiable function $F$ on $U$ such that $f = F|\sigma$.

In the classical case, $C^1[0,1] \subseteq AC(\sigma)$. We shall use the fact that all $C^1(\sigma)$ functions can be approximated by functions in $\CTPP(\sigma)$  to show that this result extends to arbitrary $\sigma$.

\section{CTPP Functions}\label{S:CTPP}

As before, we shall assume that $\sigma$ is a nonempty compact subset of the plane.
We will say that a function $f:\sigma\to\mC$ is \textbf{planar} if there exist $a,b,c\in\mC$ with $f(x,y)=ax+by+c$.
If $a,b,c\in\mR$, then (using property (\ref{const-lines}) of $BV(\sigma)$ given in Section~\ref{S:def})
\begin{equation}\label{pl-var}
	\var(f,\sigma)=\max_\sigma f- \min_\sigma f = \max_{\vecx,\vecw \in \sigma} |f(\vecx) - f(\vecw)|.
\end{equation}
This is no longer valid if $f$ is complex-valued. For complex-valued planar functions, one can split $f$ into its real and imaginary components, and find that
\begin{align*}
	\max_\sigma|f|-\min_\sigma|f|\leq\var(f,\sigma) \leq \var\big(\text{Re}(f),\sigma\big)+\var\big(\text{Im}(f),\sigma\big).
\end{align*}
Applying (\ref{pl-var}) to the real and imaginary parts shows that $f$ is of bounded variation. Note that $f$ is planar if and only if $\text{Re}(f)$ and $\text{Im}(f)$ are both planar.

We shall say that $P\subseteq\mR^2$ is a polygon if it is a compact simply connected set whose boundary consists of a finite number of line segments. By the Two Ears Theorem, all polygons can be triangulated.

Suppose that $P$ is a polygon in $\mR^2$ and let $\mathcal{A}=\{A_i\}_{i=1}^n$ be a triangulation of $P$. That is, the triangles $A_i$ are proper, closed, and have pairwise disjoint interiors, and $\cup_i A_i = P$. We say that $F:P\to \mC$ is \textbf{triangularly piecewise planar over} $\mathcal{A}$ if $F|A_i$ is planar for each $A_i$. The set of all such functions will be denoted by $\CTPP(P,\mathcal{A})$, and the set of \textbf{continuous and triangularly piecewise planar functions} on $P$ is defined to be
\begin{align*}
	\CTPP(P)=\bigcup_{\mathcal{A}}\CTPP(P,\mathcal{A}).
\end{align*}
For a given triangulation $\mathcal{A}$ of a polygon $P$, note that $f\in\CTPP(P,\mathcal{A})$ is necessarily continuous since if $A_i$ and $A_j$ are two triangles in $\mathcal{A}$, then $f|A_i$ and $f|A_j$ must agree on $A_i\cap A_j$.

If $\mathcal{A}$ and $\mathcal{A}'$ are two triangulations of $P$, we shall say that $\mathcal{A}$ is a refinement of $\mathcal{A}'$ if every triangle  in $\mathcal{A}'$ is the union of triangles in $\mathcal{A}$. Clearly in this case, if $f \in \CTPP(P,\mathcal{A}')$ then $f \in \CTPP(P,\mathcal{A})$. Any two triangulations will admit a common refinement.

\begin{definition}
	A function $f:\sigma \to \mC$ is continuous and triangularly piecewise planar if there exists a polygon $P$ containing $\sigma$, and $F \in \CTPP(P)$ such that $F|\sigma=f$. The set of all such functions will be denoted by $\CTPP(\sigma)$.
\end{definition}

If a suitable polygon $P$ exists, one can if fact choose any polygon which contains $\sigma$ as the basis for the triangulation.

\begin{lemma}\label{CTPP1}
	Suppose that $f \in \CTPP(\sigma)$. If $P_0$ is a polygon containing $\sigma$, then there exists $F_0 \in \CTPP(P_0)$ such that $f=F_0|\sigma$.
\end{lemma}

\begin{proof}
	Suppose that $P_0$ is a polygon that contains $\sigma$. By definition, there exists a polygon $P$, a triangulation $\mathcal{A}=\{A_i\}_{i=1}^n$ of $P$, and $F\in\CTPP(P,\mathcal{A})$ such that $F|\sigma=f$. If $P$ contains $P_0$, then we can take $F_0=F|P_0$. Otherwise, let $R$ be a rectangle whose interior contains both $P$ and $P_0$. Then $R\backslash\text{int}(P)$ can be triangulated, resulting in a triangulation of $R$ given by $\mathcal{A}'=\{A_i\}_{i=1}^m$. Let $P_k=\bigcup_{i=1}^k A_i$. The `ear-clipping' triangulation algorithm allows us to do this in a way that for $n+1\leq k \leq m$, the triangle $A_k$ has at least one side adjoining $P_{k-1}$, and at least one side disjoint from $P_{k-1}$ (except at the vertices).
	
	One may now extend $F$ inductively, triangle by triangle. If $F$ is already defined on $P_{k-1}$, then there exists a planar function on $A_k$ which agrees with $F$ on the intersection $A_k\cap P_{k-1}$. The triangulation $\mathcal{A}'$ generates a triangulation $\mathcal{A}_0$ of $P_0$ as follows. For each $i$, the set $A_i\cap P_0$ is a union of polygons, and thus can be written as a union of triangles. So $F_0=F|P_0\in\CTPP(P_0,\mathcal{A}_0)$ and $F_0|\sigma=f$.
\end{proof}

\begin{lemma}\label{CTPP2}
	$\textup{CTPP}(\sigma)$ is a vector space.
\end{lemma}

\begin{proof}
Closure under scalar multiplication is clear. Suppose that $f,g\in\CTPP(\sigma)$, and that $P$ is a polygon containing $\sigma$. By Lemma~\ref{CTPP1}, there exist $F,G\in\CTPP(P)$ such that $F|\sigma=f$ and $G|\sigma=g$. Clearly $F+G$ is continuous and, by using a common triangulation for $F$ and $G$, is planar on polygonal regions of $R$. Hence $F+G\in\CTPP(R)$ as required.
\end{proof}

\vfill\eject

\mysubsection{Real $\CTPP$ functions}

Some estimates are easier for real-valued functions, so we shall start by dealing with this case first, and then later give the corresponding results for complex-valued functions. We shall denote the real-valued functions in $\CTPP(\sigma)$ and $AC(\sigma)$ by $\CTPP_\mR(\sigma)$ and $AC_\mR(\sigma)$ respectively.

Suppose that $\sigma$ is a subset of a polygon $P$, that $f \in \CTPP_\mR(\sigma)$, that $F \in \CTPP_\mR(P)$ is an extension of $f$ to $P$. Since $f$ is continuous, it attains a maximum value $M$ and a minimum value $m$ on $\sigma$. The function ${\tilde F}: P \to \mR$
  \[ {\tilde F}(\vecx) = \begin{cases}
     M,   & \text{if $F(\vecx) > M$,} \\
     m,   & \text{if $F(\vecx) < m$,} \\
     F(\vecx), & \text{otherwise}
     \end{cases} \]
is also an extension of $f$ in $\CTPP_\mR(P)$ (with respect to a possibly refined triangulation). Thus a suitable extension of $f$ to $P$ can always be chosen to have the same range and the same supremum norm. For a continuous function $f:\sigma \to \mC$, let $\Delta_\sigma(f) = \sup \{|f(\vecx) - f(\vecx')| \st \vecx,\vecx' \in \sigma\}$. Obviously $\Delta_\sigma(f) \le 2 \norm{f}_\infty$.

We begin by showing that all $\CTPP_\mR$ functions are Lipschitz, and hence of bounded variation. For a triangle $A$, let $r_A$ denote the inradius of $A$. Given a triangulation $\mathcal{A} = \{A_j\}$ of a rectangle $R$, let $r(\mathcal{A}) = \min_j r_{A_j}$.

\begin{lemma}\label{CTPP3}
	Suppose that $A$ is a triangle and that $F:A\to\mR$ is planar. Then $L_A(F)$, the Lipschitz constant of $F$ on $A$, satisfies
	\[
		L_A(F)\leq\frac{\Delta_A(F)}{2 r_A}.
	\]
\end{lemma}

\begin{proof}
	Since $F$ is planar, $\nabla F$ is constant and, for distinct $\vecx,\vecx'\in A$,
	\[
		\frac{|F(\vecx)-F(\vecx')|}{\|\vecx-\vecx'\|}=\left|\nabla F\cdot \left( \frac{\vecx-\vecx'}{\|\vecx-\vecx'\|}\right)\right|.
	\]
It follows that there is a  unit vector $\vecu$ so that $L_A(F)=|\nabla F\cdot\vecu|$. Choose $\vecx,\vecx'$ on the incircle of $A$ so that $\vecx-\vecx' = 2 r_A \vecu$. Then
	\[
  L_A(F) = \frac{1}{2r_A} \left| \nabla F \cdot (\vecx - \vecx') \right|
   = \frac{1}{2r_A} | F(\vecx) - F(\vecx') |
  \leq\frac{\Delta_A(F)}{2r_A}.
	\]
\end{proof}

Note that Lemma~\ref{CTPP3} can be improved slightly, but clearly the bound must depend on the smallest width of the triangle.

\begin{theorem}\label{ctpp-lip-bv}
	$\CTPP_\mR(\sigma) \subseteq \Lip(\sigma)\subseteq BV(\sigma)$.
\end{theorem}

\begin{proof}
	Let $R$ be a rectangle that contains $\sigma$ and suppose that $f \in \CTPP_\mR(\sigma)$. By Lemma~\ref{CTPP1}, there exists a triangulation $\mathcal{A}=\{A_j\}_{j=1}^n$ of $R$ and $F\in\CTPP_\mR(R,\mathcal{A})$ such that $F|\sigma=f$. Clearly $L_R(F) \ge \max_j L_{A_j}(F)$.
	
Suppose $\vecx,\vecx'\in R$ are distinct. The line segment $\ls[\vecx,\vecx']$ joining $\vecx$ and $\vecx'$ can be written as a union of finitely many subsegments, denoted $\ls[\vecx_{j-1},\vecx_{j}]$ (with $j=1,\dots,m$), with each subsegment entirely contained in a single triangle $A_j$ in $\mathcal{A}$. Then, using Lemma~\ref{CTPP3},
	\begin{align*}
		|F(\vecx)-F(\vecx')|&	\leq\sum_{j=1}^m|F(\vecx_j)-F(\vecx_{j-1})|\\
		&	\leq\sum_{j=1}^m L_{A_j}(F) \norm{\vecx_j-\vecx_{j-1}}\\
        & \leq \max_j L_{A_j}(F) \sum_{j=1}^m \norm{\vecx_j-\vecx_{j-1}} \\
        & = \max_j L_{A_j}(F) \norm{\vecx-\vecx'}
	\end{align*}
and so $L_\sigma(f) \leq L_R(F) = \max_j L_{A_j}(F)$. Thus, $\CTPP_\mR(\sigma)\subseteq \Lip(\sigma)$. The second inclusion is just  Proposition~\ref{Lip-functs}.
\end{proof}

Note that the Lipschitz constant of $f$ may be strictly smaller than the constant for an extension $F$. It is not always possible to choose an extension $F$ with $L_\sigma(f) = L_R(F)$. 

\begin{example}\label{bad-bound} Fix $\alpha > 0$, let $\sigma = \{(x,y) \st 0 \le x \le \alpha,\ y = x^2\}$, and let $f \in \CTPP_\mR(\sigma)$ be $f(x,y) = y$. Clearly any $\CTPP_\mR$ extension $F$ to a polygon $P$  must have $L_P(F) \ge 1$. On the other hand an elementary calculation shows that $L_\sigma(f) = 2\alpha/\sqrt{1+4\alpha^2}$, which, for small $\alpha$, will be much smaller than $L_P(F)$.
\end{example}

\begin{corollary}\label{cttp-bound}
Suppose that $f \in \CTPP_\mR(\sigma)$ has an extension $F \in \CTPP_\mR(R,\mathcal{A})$ to a rectangle $R$ containing $\sigma$, and assume that $\Delta_R(F) = \Delta_\sigma(f)$.  Then
   \[ \normbv{f} \le \left(1+\frac{\ell+h}{r(\mathcal{A})}\right) \norm{f}_\infty \]
where $\ell$ and $h$ are the side lengths of $R$.
\end{corollary}

\begin{proof} Using Proposition~\ref{Lip-functs} and Lemma~\ref{CTPP3}
 \begin{align*} 
 \var(f,\sigma) \le \var(F,R) &\le C_R L_R(F) \\
   &\le (\ell + h) \max_j \frac{\Delta_{A_j}(F)}{2r_{A_j}} \\
   & \le (\ell+h) \norm{F}_\infty \max_j \frac{1}{r_{A_j}} \\
   & = \frac{(\ell+h)}{r(\mathcal{A})} \norm{f}_\infty.
 \end{align*}
\end{proof}

We also note that the bound that appears in Corollary~\ref{cttp-bound} may be far from sharp. 

\begin{example}
Let $\alpha$, $f$ and $\sigma$ be as in Example~\ref{bad-bound}.
Then $\var(f,\sigma) = \alpha^2$ and so $\norm{f}_{BV(\sigma)} = 2\alpha^2$. Let $R$ be the rectangle with corners at $(0,0),(\alpha,0),(\alpha,\alpha^2)$ and $(0,\alpha^2)$, which can be split into two triangles along the diagonal from the origin to $(\alpha,\alpha^2)$, and let $F(x,y) = y$ be the natural extension of $f$ to $R$. For this triangulation, $r(\mathcal{A}) = \frac{\alpha^2}{2 (\alpha+1)}$, so calculating the bound above
gives 
  \[ \left(1+ \frac{\ell+h}{r(\mathcal{A})} \right) \norm{f}_\infty
   = \left(1+ \frac{2(\alpha+1)(\alpha+\alpha^2)}{\alpha^2}\right) \alpha^2
     = 2 \alpha + 5 \alpha^2 + 2\alpha^3 \]
which will be much larger than $\norm{f}_{BV(\sigma)}$ if $\alpha$ is small.
\end{example}



\section{$\CTPP_\mR(\sigma)$ is contained in $AC_\mR(\sigma)$}

Our next step is to show that every function in $\CTPP_\mR(\sigma)$ is not just of bounded variation, but is absolutely continuous.
Suppose that $f \in \CTPP_\mR(\sigma)$ and that $F$ is a piecewise planar extension of $f$ to some rectangle $R$ containing $\sigma$. If $F$ is absolutely continuous, then so is its restriction $f$. Therefore it will be sufficient to just consider $\CTPP$ functions defined on a rectangle.

Suppose then that $R$ is a rectangle in the plane and that $f \in \CTPP_\mR(R)$ with respect to a fixed triangulation $\mathcal{A}$ of $R$. Our aim is to show that each point $\vecx \in R$ admits a compact neighbourhood $U_\vecx$ on which $f|U_\vecx$ is absolutely continuous and then to use the Patching Lemma (Theorem~\ref{patching-lemma}).

Note that since we can choose $R$ as large as we like, it is sufficient to do this for each interior point of $R$ as this will ensure that $f$ is absolutely continuous on any compact subset of the interior of $R$. Doing this avoids some special cases in the argument below.


\begin{definition}
We will say that an interior point $\vecx \in R$ is
\begin{itemize}
	\item a \textbf{planar point} for $f$ if it lies in exactly one triangle of $\mathcal{A}$;
	\item an \textbf{edge point} for $f$ if it lies in exactly two triangles of $\mathcal{A}$;
	\item a \textbf{vertex point} for $f$ if it lies in three or more triangles of $\mathcal{A}$.
\end{itemize}
\end{definition}

Note that these three cases are exhaustive and mutually exclusive, and that the classification of $\vecx$ depends on the triangulation. Dealing with the first class of points is straightforward.


Suppose first that $\vecx$ is a planar point lying in the triangle $A \in \mathcal{A}$. In this case $\vecx$ is an interior point of $A$ and so  $U_\vecx =  A$ is a compact neighbourhood of $\vecx$. As $f$ is planar on $U_\vecx $, it is absolutely continuous on this compact neighbourhood of $\vecx$.


\begin{lemma}\label{edge-pt}
	Suppose that $\vecx$ is an edge point for $f \in \CTPP_\mR(R)$ with respect to the triangulation $\mathcal{A}$. Then there exists a compact neighbourhood $U_{\vecx}$ of $\vecx$ such that $f|U_{\vecx} \in AC(U_\vecx)$.
\end{lemma}

\begin{proof}
	By affine invariance, we may assume that $\vecx=(0,0)$. Moreover, we may assume that the two triangles that $\vecx$ lies in intersect on the line $x=0$. Thus, there exists some small square $R_0$ centred at the origin such that if $(x,y)\in R \cap R_0$, then
	\begin{align*}
		f(x,y)=\begin{cases}
			a_1x+b_1y+c_1 & \text{if $x \geq 0$},\\
			a_2x+b_2y+c_2 & \text{if $x \leq 0$}.
		\end{cases}
	\end{align*}
Since $f$ is continuous, $b_1=b_2$ and $c_1=c_2$. Let $U_\vecx=R \cap R_0$ and
	\begin{equation*}
		f_1(x,y)=\begin{cases}
			a_1x & \text{if $x \geq 0$},\\
			a_2x & \text{if $ x \leq 0$},
		\end{cases}
		\qquad \qquad f_2(x,y)=b_1y+c_1
	\end{equation*}
so that $f = f_1+f_2$.
It follows from Property~(\ref{const-lines}) in Section~\ref{S:def} that $f_1,f_2\in AC(U_\vecx)$, and so $f\in AC(U_\vecx)$.	
\end{proof}

The most difficult case is that of vertex points. Let $Q=Q_t=[-t,t]^2\subset\mR^2$. We shall say that a function $f:Q \to \mR$ is \textbf{star-planar} on $Q$ if there exists a partitioning $\mathcal{A}=\{A_i\}_{i=1}^n$ of $Q$ given by $n\geq2$ rays starting at the origin such that $f$ is planar on each set $A_i$. (If desired, one may enforce that each $A_i$ is triangular by adding the rays from the origin to the four corners of $Q$.)

\begin{lemma}\label{star-planar}
	If $f$ is star-planar on $Q$ with respect to $\{A_i\}_{i=1}^n$, then $f\in BV(Q)$ with
	\begin{align*}
		\var(f,Q)\leq 2n\sup_{\vecx,\vecw\in Q}|f(\vecx)-f(\vecw)|.
	\end{align*}
\end{lemma}

\begin{proof}
	First note that if $f$ is star-planar on $Q$ with respect to $\{A_i\}_{i=1}^n$, then it is star-planar with respect to a finer partition $\{A_i'\}_{i=1}^m$ with $m\leq 2n$ formed by extending the rays to full lines. The illustration below provides an example of such an extension.
	\begin{center}
		\begin{tikzpicture}
			\draw (0,0)--(0,2)--(2,2)--(2,0)--(0,0);
			\draw (3,1)--(5,1)--(4.8,1.2)--(5,1)--(4.8,0.8);
			\draw (6,0)--(6,2)--(8,2)--(8,0)--(6,0);
			
			\draw (1,1)--(0,0.5);
			\draw (1,1)--(0.8,2);
			\draw (1,1)--(2,1.1);
			
			\draw (8,1.5)--(6,0.5);
			\draw (7.2,0)--(6.8,2);
			\draw (6,0.9)--(8,1.1);
		\end{tikzpicture}
	\end{center}
Without loss, assume that $A_i'$ shares an edge with $A_{i+1}'$ for $1 \le i< m$.
As $f$ is planar on $A_i'$, Equation (\ref{pl-var}) gives that
	\begin{align*}
		\var(f,A_i)=\sup_{\vecx,\vecw\in A_i}|f(\vecx)-f(\vecw)|.
	\end{align*}
The finer partition and ordering of $A_i'$s has the property that $B_k = \bigcup_{i=1}^{k}A_i'$ is convex for each $1\leq k\leq \frac{m}{2}$. By using Theorem~\ref{join-convexly} inductively, one deduces that $f|_{B_k}\in BV(B_k)$ for each $1\leq k\leq \frac{m}{2}$. Similarly, letting $C_k = \bigcup_{i=0}^{k-1} A_{m-i}'$ for $1\leq k\leq \frac{m}{2}$, we find that $f|_{C_k}\in BV(C_k)$. Since $Q = B_{m/2} \cup C_{m/2}$ is convex, using Theorem~\ref{join-convexly} one last time gives $f\in BV(Q)$ with the claimed bound.
\end{proof}


\begin{lemma}\label{vertex-pt}
	Suppose that $\vecx$ is a vertex point for $f \in \CTPP_\mR(R)$ with respect to the triangulation $\mathcal{A}$. Then there exists a compact neighbourhood $U_\vecx$ of $\vecx$ such that $f|U_\vecx \in AC(U_\vecx)$.
\end{lemma}

\begin{proof}
	By affine invariance, we may assume that $\vecx=(0,0)$. Also, by replacing $f$ with $f-f(0,0)$, we may assume that $f(0,0)=0$.
	

As above, let $Q_t$ be the square $[-t,t]^2$. One may choose $\delta>0$ so that $Q_\delta$ lies entirely within $R$ and contains no other vertex points for $f$. Thus $f$ is star-planar on $Q_\delta$ with respect to the partitioning induced by $\mathcal{A}$. Let $U_\vecx=Q_\delta$. To show that $f|U_\vecx$ can be approximated by absolutely continuous functions, it suffices to show that for all $\varepsilon>0$, there is an $h \in AC(Q_\delta)$ with $\norm{f-h}_{BV(Q_\delta)} < \varepsilon$.
	
Fix $\varepsilon > 0$. By the uniform continuity of $f$ and Lemma~\ref{star-planar}, one may find $s \in (0,\delta)$ so that $\norm{f}_{BV(Q_s)} < \varepsilon/10$. Define $g_s: [-\delta,\delta] \to \mR$ as the piecewise linear function whose graph is shown below.
	\begin{center}
		\begin{tikzpicture}[scale=0.7]
			\draw (-6.5,0)--(6.5,0);
			\draw (0,-1)--(0,2.5);
			
			\draw (-6,0.2)--(-6,-0.2) node[below] {$-\delta$};
			\draw (6,0.2)--(6,-0.2) node[below] {$\delta$};
			
			\draw (4,0.2)--(4,-0.2) node[below] {$s$};
			\draw (-4,0.2)--(-4,-0.2) node[below] {$-s$};
			
			\draw (-2,0.2)--(-2,-0.2) node[below] {$-s/2$};
			\draw (2,0.2)--(2,-0.2) node[below] {$s/2$};
			
			\draw[blue,very thick] (-6,2)--(-4,2)--(-2,0)--(2,0)--(4,2)--(6,2) node[black, below right] {$g_s$};
			
			\draw (-0.2,2)--(0.2,2) node[right] {$1$};
		\end{tikzpicture}
	\end{center}
Note that $g_s\in AC[-\delta,\delta]$ with $\norm{g_s}_{BV[-\delta,\delta]} = 3$. By Property~(\ref{const-lines}), the functions $(x,y)\mapsto g_s(x)$ and $(x,y)\mapsto g_s(y)$ are both in $AC(Q_\delta)$ with $BV(Q_\delta)$ norm equal to three. It follows that their product $\chi$ is also in $AC(Q_\delta)$ with $\norm{\chi}_{BV(Q_\delta)}\leq10$.
	
	Let $h=\chi f$, let $B = Q_\delta \backslash (-s/3,s/3)^2$, and suppose that $\vecw \in B$. As $Q_\delta$ contained only one vertex point, namely the origin, then $\vecw$ is either a planar point or an edge point for $f$. As shown earlier, there is a compact neighbourhood $V_\vecw$ of $\vecw$ such that $f|V_\vecw\in AC(V_\vecw)$. Clearly $\chi|V_\vecw\in AC(V_\vecw)$, and so $h|V_\vecw\in AC(V_\vecw)$. On the other hand, if $\vecw\in (-s/3,s/3)^2$, then $h=0$ on an open neighbourhood of $\vecw$, and so again we can choose a compact neighbourhood $V_\vecw$ of $\vecw$ such that $h|V_\vecw\in AC(V_\vecw)$. It follows from the Patching Lemma (Theorem~\ref{patching-lemma}) that $h\in AC(Q_\delta)$. Moreover, we have that
	\[
   \norm{f-h}_{BV(Q_\delta)}
     = \norm{f(1-\chi)}_{BV(Q_{\delta})}
     \leq \norm{f}_{BV(Q_\delta)} \norm{1-\chi}_{BV(Q_\delta)}
     < \varepsilon.
	\]
Thus, $f\in AC(Q_\delta)$.
\end{proof}

Combining these results with the Patching Lemma gives the following corollary.

\begin{corollary}\label{CTPP-AC}
	$\CTPP_\mR(\sigma)\subseteq AC_\mR(\sigma)$.
\end{corollary}

\section{$C^1$  functions are absolutely continuous}\label{S:C1}

\begin{theorem}\label{CTPP-dense1}
	Let $R$ be a closed rectangle in the plane. Then $C^1_\mR(R)\subseteq \textup{cl}(\CTPP_\mR(R))$, and hence $C^1_\mR(R)\subseteq AC_\mR(R)$.
\end{theorem}

\begin{proof}
Note that both the $C^1$ functions and the absolutely continuous functions are preserved under affine transformations of the plane, so it suffices to prove the theorem for $R = [0,1]^2$.

Suppose then that $R = [0,1]^2$ and that $f \in C^1_\mR(R)$. Fix a $C^1$ extension of $f$ to an open neighbourhood of $R$.
Fix $\varepsilon > 0$ and (by the uniform continuity of $f$ and $\nabla f$) choose $\delta > 0$ such that for all $\vecx,\vecy \in R$ with $\norm{\vecx-\vecy} \le \delta$,
	\begin{equation}\label{bounds}
		|f(\vecx)-f(\vecy)| < \varepsilon
        \qquad \text{and} \qquad
         \norm{\nabla f(\vecx)-\nabla f(\vecy)} < \varepsilon.
	\end{equation}
Let $n=\left\lceil \sqrt{2}/\delta \right\rceil$ and let $\mathcal{A}$ be the triangulation of $R$ shown in the figure below.
	\begin{center}
		\begin{tikzpicture}
			\draw[blue] (0,0)--(4,0)--(4,4)--(0,4)--(0,0);
			\draw (0,0)--(4,4);
			\draw (1,0)--(4,3);
			\draw (2,0)--(4,2);
			\draw (3,0)--(4,1);
			\draw (0,1)--(3,4);
			\draw (0,2)--(2,4);
			\draw (0,3)--(1,4);
			
			\node[left] at (0,4) {$1$};
			\node[below] at (4,0) {$1$};
			\draw (4,3)--(0,3) node[left=3] {$\vdots$};
			\draw (4,2)--(0,2) node[left] {$\frac{2}{n}$};
			\draw (4,1)--(0,1) node[left] {$\frac{1}{n}$};
			\draw (3,4)--(3,0) node[below=6] {$\hdots$};
			\draw (2,4)--(2,0) node[below] {$\frac{2}{n}$};
			\draw (1,4)--(1,0) node[below] {$\frac{1}{n}$};
		\end{tikzpicture}
	\end{center}
Note that each triangle has diameter less than or equal to $\delta$. Let $g\in\CTPP(R,\mathcal{A})$ agree with $f$ at all vertices in the above triangulation.
	
Fix a triangle $A \in \mathcal{A}$. The bounds in (\ref{bounds}) imply that there exist $m,M$ such that $M-m<\varepsilon$ and $m\leq f\leq M$ on $A$. As $g$ is planar on $A$ (so its extrema occur at the vertices of $A$) and $g$ agrees with $f$ at the vertices of $A$, we also have that $m\leq g\leq M$. So for any $\vecx\in A$,
	\begin{align*}
		|f(\vecx)-g(\vecx)| \leq |f(\vecx)-m|+|g(\vecx)-m|
                < 2 \varepsilon.
	\end{align*}
	Thus $\norm{f-g}_\infty < 2 \varepsilon$.
	
	We will now estimate the Lipschitz constant of $d=f-g$. Suppose that $\vecx,\vecy\in A$ are distinct, and let $\ell$ denote the line segment from $\vecx$ to $\vecy$. By the Mean Value Theorem, there exists $\vecq \in \ell$ such that
	\[
		|d(\vecx)-d(\vecy)| \leq \norm{\nabla d(\vecq)}\, \norm{\vecx-\vecy}.
	\]
The planarity of $g$ implies that $\nabla g$ is constant. Using the fact that $f=g$ on the vertices of $A$, applying Rolle's Theorem to $d$ along the sides of $A$ parallel to the coordinate axes yields $\boldsymbol\xi$ and $\boldsymbol\eta$ on the boundary of $A$ such that
	\[
		\nabla g(\vecq) = \left( \frac{\partial f}{\partial x}(\boldsymbol\xi) ,
                  \frac{\partial f}{\partial y}(\boldsymbol\eta) \right).
	\]
Thus we can estimate $\|\nabla d(\vecq)\|$ as follows:
	\[
	\norm{\nabla d(\vecq)}^2
       = \left( \frac{\partial f}{\partial x}(\vecq)-\frac{\partial f}{\partial x}(\boldsymbol\xi) \right)^2
       + \left( \frac{\partial f}{\partial y}(\vecq)-\frac{\partial f}{\partial y}(\boldsymbol\eta) \right)^2
       < 2 \varepsilon^2
	\]
	by (\ref{bounds}), and so we have that
	\[
		|d(\vecx)-d(\vecy)| < \sqrt{2}\,\varepsilon \norm{\vecx-\vecy}.
	\]
	One may extend this last inequality to general distinct $\vecx,\vecy\in R$ by splitting the line segment between $\vecx$ and $\vecy$ into segments which lie entirely in a single triangle in the spirit of the proof of Theorem~\ref{ctpp-lip-bv}. Thus, using Proposition~\ref{Lip-functs},
	\[
		\norm{f-g}_{BV(R)} \le C_R \norm{f-g}_{\Lip(R)}
            < 2 (2+\sqrt{2}) \varepsilon.
	\]
Thus $C^1_\mR(R)\subseteq \textup{cl}(\CTPP_\mR(R))\subseteq AC_\mR(R)$.
\end{proof}

\begin{theorem}\label{CTPP-dense3}
	Suppose that $\sigma$ is a nonempty compact subset of the plane. Then $C^1_\mR(\sigma)$ is a dense subset of $AC_\mR(\sigma)$.
\end{theorem}

\begin{proof}
	Suppose that $f\in C^1(\sigma)$, and thus admits a $C^1$ extension (also denoted by $f$) on some open neighbourhood $V$ of $\sigma$. Suppose that $\vecx\in\sigma$. Then there exists a closed rectangle $R_\vecx$ centred at $\vecx$ which lies inside $V$. By Theorem~\ref{CTPP-dense1}, $f|R_\vecx\in AC(R_\vecx)$. The set $U_\vecx=R_\vecx\cap\sigma$ is a compact neighbourhood of $\vecx$ such that $f|U_\vecx\in AC(U_\vecx)$. By the Patching Lemma, one deduces that $f\in AC(\sigma)$. The density of $C^1(\sigma)$ follows from the fact that the polynomials are contained in $C^1(\sigma)$.
\end{proof}

\begin{theorem}\label{CTPP-dense2}
	$\CTPP(\sigma)$ is dense in $AC(\sigma)$.
\end{theorem}

\begin{proof}
	Suppose that $f\in AC(\sigma)$. Fix a closed rectangle $R$ containing $\sigma$. For a given $\varepsilon>0$, there exists a polynomial $p$ such that $\norm{f-p}_{BV(\sigma)} < \varepsilon/2$. As $p \in C^1(R)$, Theorem~\ref{CTPP-dense1} implies the existence of $g \in \CTPP(R)$ such that $\norm{p-g}_{BV(R)} < \varepsilon/2$. Hence $\norm{f-g}_{BV(\sigma)} < \varepsilon$.
\end{proof}

\section{Complex-valued functions}

Suppose now that $f \in AC(\sigma)$ is complex-valued. As noted earlier, if $f_R = \Re(f)$ and $f_I = \Im(f)$, then $f_R,f_I \in \AC_\mR(\sigma)$. Suppose then that $\epsilon > 0$. By Theorem~\ref{CTPP-dense2}, there exist $g_R,g_I \in \CTPP_\mR(\sigma)$ with $\normbv{f_R-g_R} < \frac{\epsilon}{2}$ and $\normbv{f_I-g_I} < \frac{\epsilon}{2}$. By Theorem~\ref{CTPP2}, $g = g_R+ig_I \in \CTPP(\sigma)$, and $\normbv{f-g} < \epsilon$.

The same argument shows that there exists $h \in C^1(\sigma)$ with $\normbv{f-h} < \epsilon$.

\begin{theorem} $\CTPP(\sigma)$ and $C^1(\sigma)$ are dense in $AC(\sigma)$.
\end{theorem}

\section*{Acknowledgements}
The work of the third author was supported by the Research Training Program of the Department of Education and Training of the Australian Government.

%
%
\bibliographystyle{amsalpha}

\end{document}